\documentclass[11pt]{amsart}
\usepackage{amssymb}
\usepackage{mathrsfs}
\usepackage[english]{babel}
\swapnumbers

\begin{document}
\newtheorem*{thm*}{Theorem}
\newtheorem*{cor*}{Corollary}
\def\id{{\it I\! d}}
\def\bbb#1{\overline{\!#1}}
\def\se{\subseteq}
\def\d{\,\mathrm{d}}
\def\CC{\mathbf{C}}
\def\HH{\mathbf{H}}
\def\RR{\mathbf{R}}
\def\ZZ{\mathbf{Z}}
\def\F{\mathscr{F}}
\def\O{\mathscr{O}}
\def\ro{\varrho}
\def\fhi{\varphi}
\def\teta{\vartheta}
\def\epsi{\varepsilon}
\def\wt{\widetilde}
\def\ti{-\allowhyphens}
\def\lra{\longrightarrow}
\def\No{N\raise4pt\hbox{\tiny o}\kern+.2em}
\def\no{n\raise4pt\hbox{\tiny o}\kern+.2em}
\def\bsl{\backslash}
\def\beq{\begin{equation}}
\def\eeq{\end{equation}}
%


\title[Arithmeticity \emph{vs.} Non-Linearity]{Arithmeticity \emph{vs.} Non-Linearity \\ for Irreducible Lattices}
\author[Nicolas Monod]{Nicolas Monod*}
\thanks{*University of Chicago, {\sf monod@math.uchicago.edu}}
\begin{abstract}
We establish an arithmeticity \emph{vs.} non-linearity alternative for irreducible lattices in suitable product groups, for instance products of topologically simple groups. This applies notably to a (large class of) Kac-Moody groups. The alternative relies heavily on the superrigidity theorem we propose in~\cite{MonodCAT0}, since we follow Margulis' reduction of arithmeticity to superrigidity.
\end{abstract}
\maketitle

\section{Introduction}

\subsection{Arithmeticity}
Margulis proved that every irreducible lattice in a Lie/algebraic group of higher rank is \emph{arithmetic}; see~\cite[IX]{Margulis} for exact statements (some authors write \emph{$S$\ti arithmetic} in the generality that we consider). The main step consists in establishing his \emph{superrigidity} theorems for linear representations.

We presented in~\cite{MonodCAT0} superrigidity results for a wide class of lattices in quite arbitrary locally compact product groups. These results apply also, as a rather particular case, to linear representations. Therefore, we can deduce certain general arithmeticity statements from our superrigidity following Margulis' line of argumentation; this is the object of the present note. That goal requires some assumptions on the ambient groups; topological simplicity is more than sufficient here, and we shall propose a wider framework, see~\ref{sec_few_factors}.

\medskip

\subsection{Non-Linearity}
The context of lattices in general product groups is particularly interesting because there are remarkable examples outside the arithmetic world; for instance, the groups of Burger-Mozes~\cite{Burger-Mozes1},\cite{Burger-Mozes2} and Kac-Moody groups as studied by R\'emy~\cite{RemyAST},\cite{RemyERN},\cite{Remy-Bonvin},\cite{RemySURV} and Carbone-Garland~\cite{Carbone-Garland}. Several of these groups are known to be non-linear~\cite{Burger-Mozes2},\cite{Remy-Bonvin}; and indeed, reversing the viewpoint on our arithmeticity results below, it follows for suitable groups that unless the lattice is arithmetic and the topological groups algebraic, \itshape the lattice does not admit any linear representation over any field (of characteristic~$\neq 2,3$) with infinite image.\upshape

\medskip

\subsection{Comments}
(a)~Shalom's cohomological superrigidity~\cite{Shalom00} allowed him notably to deal with the case of representations to $\mathbf{GL}_2$, and deduce a corresponding arithmeticity statement. He complemented Margulis' strategy by a result of Pink~\cite{Pink98} and his superrigidity for characters~\cite[0.8]{Shalom00}; we shall borrow these two arguments aswell.

\smallskip

(b)~In the particular case of lattices in products of tree groups, Burger-Mozes-Zimmer have obtained detailed superrigidity/arithmeticity results in characteristic zero (in~\cite{Burger-Mozes-Zimmer}; see announcement in~\cite{MozesICM}).

\smallskip

(c)~A recent result of Bader-Shalom~\cite{Bader-Shalom} establishes a ``normal subgroup theorem'' \emph{\`a la} Margulis for lattices as those under consideration here; we will however not need this result in our proof.

\smallskip

(d)~Even though we also provide commensurator superrigidity in~\cite{MonodCAT0}, we do not consider that case below as it has already been treated quite extensively by other authors (\emph{e.g.}~\cite{Margulis},\cite{MargulisCOM},\cite{A'Campo-Burger},\cite{BurgerICM},\cite{Remy-Bonvin},\cite{Shalom00}).

\medskip

\subsection{Acknowledgments}
As already mentioned, the reduction of arithmeticity to superrigidity follows closely Margulis' argumentation~\cite[IX]{Margulis} and owes also to Shalom~\cite[\S8]{Shalom00}, especially for~\ref{sec_large_NA}. I am grateful to D.~Fisher and Y.~Shalom for encouraging me to write this note and for helpful comments; to B.~Klingler and H.~Oh for helpful comments; to D.~Witte Morris for the reference to~\cite{Hochschild81}.

This work was partially supported by FNS grant~8220-067641 and NSF grant DMS~0204601.

\medskip

\section{Definitions and Statements of Results}

\subsection{}
\label{sec_G_i}%
We shall consider throughout this paper locally compact groups
$$G = G_1\times \cdots \times G_n,$$
where $n\geq 2$ and each $G_i$ is compactly generated. A \emph{lattice} $\Gamma<G$ is a discrete subgroup of finite invariant covolume; any discrete cocompact subgroup is a lattice. We call a lattice $\Gamma$ \emph{completely irreducible} if (i)~its projection to any factor $G_i$ is dense and (ii)~its projection to any (proper) subproduct of $G_i$'s is non-discrete. For lattices in semisimple Lie/algebraic groups without compact factors (where $n\geq 2$ means we are in the non-almost-simple case), these conditions are essentially equivalent to the algebraic notion of irreducibility~\cite[II.6.7]{Margulis}.

\medskip

\noindent
\textbf{Remarks.~}(a)~Condition~(i) by itself is not really a restriction; indeed, if $\Gamma$ is \emph{any} lattice in $G$, then one verifies that it still is a lattice in the product of all $G'_i$, where $G'_i$ is the closure of the projection of $\Gamma$ to $G_i$ (for the finiteness of covolume, see \emph{e.g.}~\cite[1.6]{Raghunathan}).

\smallskip

\noindent
(b)~Whilst condition~(ii) was not imposed in~\cite{MonodCAT0}, it is necessary here in order to rule out direct product lattices. For instance, one can obtain trivial counter-examples to Theorem~\ref{sec_few_factors} below by taking $\Gamma = G = G_1 \times G_2$, where $G_1$, $G_2$ are suitable infinite finitely generated groups with $G_1$ linear and $G_2$ non-linear.

\smallskip

\noindent
(c)~It follows from~(ii) that each $G_i$ is non-discrete and non-compact.

\medskip

\subsection{}
\label{sec_non_uniform}%
If $\Gamma$ is a cocompact lattice in a compactly generated locally compact group $G$, then:

\smallskip

\noindent
(i)~$\Gamma$ is \emph{square-integrable}, that is: It is finitely generated, and for the length function $\ell$ associated to some (or equivalently any) finite generating set, there is a Borel right $\Gamma$\ti fundamental domain $\F\se G$ (with null boundary) such that
$$\int_\F \ell(\chi(g^{-1} h))^2\d h \ <\infty\kern1cm\forall\,g\in G,$$
wherein $\chi:G\to\Gamma$ is uniquely defined by: $\forall\,g\in G, g\in\F\chi(g)^{-1}$.

\smallskip

\noindent
(ii)~$\Gamma$ is \emph{weakly cocompact}, that is: The $G$\ti representation $L^2_0(G/\Gamma)$ (\emph{i.e.} the orthogonal complement of the trivial representation in $L^2(G/\Gamma)$) does not almost have non-zero invariant vectors.

\smallskip

For more details, see~\cite[App.~B]{MonodCAT0}; the point is that the results in~\cite{MonodCAT0} also hold for non-cocompact lattices, provided they satisfy~(i) and~(ii) above. Thus, square-integrability and weak cocompactness are going to be assumptions in our statements below.

Besides the cocompact case,~(i) and~(ii) hold notably for many of R\'emy's Kac-Moody lattices; indeed,~(i) was recently established by R\'emy~\cite{Remy04}, and in any case~(ii) trivially holds whenever one considers groups with Kazhdan's property~(T); it follows from~\cite{Dymara-J02} that every Kac-Moody group over ${\bf F}_q$ whose Cartan matrix has finite entries is Kazhdan whenever~$q$ is large enough. For the classical Lie/algebraic case, compare~\ref{sec_classical}.

\medskip

\subsection{}
Let $K$ be a global field, $\HH$ a connected semisimple $K$\ti group, $S$ a set of inequivalent valuations of $K$; if the characteristic of $K$ is positive, then we assume $S\neq \varnothing$. A subgroup of $\HH$ is called \emph{arithmetic} if it is commensurable with $\HH(K(S))$, where $K(S)$ is the ring of $S$\ti integers of $K$; we recall that two subgroups are called \emph{commensurable} if their intersection is of finite index in both.

\medskip

The following well-known construction presents certain arithmetic groups as lattices in compactly generated locally compact groups. For any  valuation $v$, denote the corresponding completion of $K$ by $K_v$ and endow $\HH(K_v)$ with the corresponding locally compact topology. Assume that $S$ is finite and contains all Archimedean valuations $v$ for which $\HH$ is $K_v$\ti isotropic. Reduction theory (after Borel-Harish-Chandra and Behr-Harder) shows that the diagonal embedding realizes $\HH(K(S))$ as a lattice in the product $\prod_{v\in S} \HH(K_v)$. This will be our archetypal arithmetic lattice. We observe that infinite sets $S$ also give rise to lattices, but in ad\'elic groups, precluding compact generation. For all the above, we refer to~\cite[I.3]{Margulis}.

\medskip

We shall consider below the normal subgroup $\HH(K_v)^+$ of $\HH(K_v)$ generated by the $K_v$\ti points of the unipotent radical of all parabolic $K_v$\ti subgroups of $\HH$. This group is introduced in detail by Borel-Tits~\cite[\S6]{Borel-Tits73}. We shall recall its properties as we need them; for now we just recall that it is a closed cocompact normal subgroup with Abelian quotient~\cite[6.14]{Borel-Tits73} (assuming that $\HH$ has no non-trivial $K_v$\ti anisotropic factors).

\medskip

\subsection{}
\label{sec_thm_simple}%
We now give the first form of our arithmeticity \emph{vs.} non-linearity alternative. One obtains a clear-cut statement by assuming that each $G_i$ is \emph{topologically simple}, \emph{i.e.} the only closed normal subgroups of $G_i$ are $G_i$ and $\{e\}$.

\begin{thm*}
Let $G=G_1 \times \cdots\times G_n$ be a product of topologically simple compactly generated locally compact groups and $\Gamma<G$ a completely irreducible square-summable weakly cocompact lattice. Then either

\begin{itemize}
\item[(i)] There is a topological isomorphism $G\cong \prod_{v\in S} \HH(K_v)^+$ under which $\Gamma$ is commensurable with $\HH(K(S))$, where $K$ is a global field, $\HH$ a connected absolutely simple adjoint $K$\ti group, $S$ a finite set of inequivalent valuations; or:

\item[(ii)] Any homomorphism from $\Gamma$ to any linear group over any field of characteristic~$\neq2,3$ has finite image.
\end{itemize}
\end{thm*}

\medskip

\noindent
\textbf{Remarks.~}(a)~In alternative~(i), we view $\HH(K(S))$ diagonally embedded in $\prod_{v\in S} \HH(K_v)$, so commensurability with the image of $\Gamma$ in $\prod_{v\in S} \HH(K_v)^+$ makes sense \emph{via} the inclusions $\HH(K_v)^+ \to \HH(K_v)$. In addition, see~\ref{sec_classical}.

\smallskip

\noindent
(b)~In alternative~(i), $S$ contains all the Archimedean valuations $v$ for which $\HH$ is $K_v$\ti isotropic.

\smallskip

\noindent
(c)~The assumption on the characteristic is used in several points in the argument to rule out algebraic groups with non-standard isogenies and their pathologies; we do not expect this assumption to be indispensable, see~\cite{Venkataramana88} for the classical case.

\medskip

The various assumptions on $\Gamma$ in the above theorem might have obscured the simplicity of the statement; in order to dispel this impression, we point out the following elementary formulation (for all proofs, see Section~\ref{sec_proofs}):

\begin{cor*}
Let $\Gamma$ be any discrete cocompact subgroup of a product $G=G_1 \times G_2$ of compactly generated locally compact topologically simple groups $G_i$. If the projections of $\Gamma$ to both $G_i$ are dense, then the alternative (i)/(ii) of Theorem~\ref{sec_thm_simple} holds.
\end{cor*}

\medskip

\subsection{}
\label{sec_classical}%
At first sight, Theorem~\ref{sec_thm_simple} does not apply to the ``classical'' case of irreducible lattices in semisimple groups, since the latter need not be topologically simple. It does however apply after the following reduction steps:

First, one may immediately assume that there are no compact normal subgroups by passing to the adjoint absolutely simple case (and factoring out possible anisotropic factors). Then, according to Tits~\cite{Tits64}, we are reduced to simple groups upon replacing each factor $\HH(K_v)$ with $\HH(K_v)^+$. The only problem is that in positive characteristic the group $\HH(K_v)^+$ need not have finite index in $\HH(K_v)$. However, the quotient is Abelian and torsion~\cite[6.14]{Borel-Tits73}; therefore, the image of any finitely generated group in $\HH(K_v)/\HH(K_v)^+$ is finite and thus the lattice lies indeed in the product of the $\HH(K_v)^+$'s up to commensurability.

\smallskip

There is however one limitation of our method for non-cocompact lattices, arising from the assumption of square-integrability: Assume that $\Gamma$ is a non-cocompact irreducible lattice in a (non-almost-simple) semisimple Lie/algebraic group. The condition~(ii) of~\ref{sec_non_uniform} is taken care of by property~(T) unless there are rank one factors; moreover, in the Lie group case,~(ii) always holds even without property~(T) (see~\cite[II.1.12]{Margulis} and~\cite{Bekka98}). Condition~(i), however, is known to hold~\cite[\S2]{Shalom00} but only as an application of Margulis' arithmeticity.

\medskip

\subsection{}
\label{sec_few_factors}%
We now proceed to relax the simplicity assumption on $G_i$. A first class of groups that we can encompass is the topological analogue of \emph{hereditarily just infinite} groups, namely: 
$$\text{Every non-trivial closed normal subgroup of $G_i$, and of}\atop 
\text{its finite index open subgroups, has finite index.}\leqno(*)$$
This setting is much more general than topological simplicity, since a group $G_i$ as above can for instance be residually finite. We can relax this assumption further, though it will sound more cumbersome: We shall say that $G_i$ has \emph{few factors} if the following hold:

\begin{itemize}
\item[(i)] Every non-trivial closed normal subgroup of $G_i$ is cocompact.

\item[(ii)] There is no non-zero continuous homomorphism $G_i\to\RR$.

\item[(iii)]Every closed normal cocompact subgroup of $G_i$ satisfies~(i) and~(ii).
\end{itemize}

\noindent
Property~$(*)$ implies that $G_i$ has few factors (since it is non-discrete, see~\ref{sec_proof_few} for proof). In fact, we shall need less in our proofs, as we will only apply~(iii) to \emph{finite index} subgroups and to the identity component $G_i^0$; thus all the following holds in this more relaxed setting.

\begin{thm*}
Let $G=G_1 \times \cdots\times G_n$ be a product of compactly generated locally compact groups with few factors and $\Gamma<G$ a completely irreducible square-summable weakly cocompact lattice. Then either

\begin{itemize}
\item[(i)] After replacing $\Gamma$ and each $G_i$ with finite index (closed) subgroups, there is a topological isomorphism $G\cong \prod_{v\in S} \HH(K_v)^+$ as in Theorem~\ref{sec_thm_simple}, under which the image of $\Gamma$ is commensurable with $\HH(K(S))$; or:

\item[(ii)] Any homomorphism from $\Gamma$ to any linear group over any field of characteristic~$\neq2,3$ has finite image.
\end{itemize}
\end{thm*}

As a particular case, one deduces that there are no irreducible lattices in ``mixed'' products of some Lie/algebraic factors with some ``exotic'' (non-algebraic) groups (with few factors). We do not know whether there is a simpler proof of this fact; in any event, let us state the following particular case for the record:

\begin{cor*}
Consider a lattice $\Gamma<G = G_1\times \cdots \times G_n$ as in the theorems. If any of the $G_i$ is connected, then we are in the arithmetic alternative~(i).
\end{cor*}

Of course, the same corollary holds if any of the $G_i$ is linear (in characteristic~$\neq 2,3$), since the image of $\Gamma$ in $G_i$ is infinite by Remark~\ref{sec_G_i}(c) and density.

If $G_i$ is connected, then it follows from the solution to Hilbert's fifth problem~\cite[4.6]{Montgomery-Zippin} and from the absence of compact normal subgroups that it is a Lie group; the assumptions on $G_i$ imply that this Lie group is linear.

\section{Proofs}
\label{sec_proofs}%

\subsection{}
We shall adopt the setting of Theorem~\ref{sec_few_factors}, since it implies immediately Theorem~\ref{sec_thm_simple}; the corollary of~\ref{sec_thm_simple} is proved in~\ref{sec_cor}. Thus, we suppose that for some field $F$ of characteristic~$\neq 2,3$ and some $d\in\mathbf{N}$ there is a homomorphism $f:\Gamma\to\mathbf{GL}_d(F)$ with infinite image and proceed to show that the assertion~(i) of Theorem~\ref{sec_few_factors} holds. Observe that if we replace $\Gamma$ with a finite index subgroup and each $G_i$ with the closure of the projection of our finite index subgroup to $G_i$, then the latter is of finite index in $G_i$ (and hence open); one may further reduce to the case where $G_i$ is normal. Moreover, one verifies that after this operation all other assumptions are preserved, both regarding $G$ and the lattice. The conclusion of Theorem~\ref{sec_few_factors} is also stable under this operation; we will perform this below in~\ref{sec_start} and~\ref{sec_reduce_Gamma}.

\smallskip

We shall attempt to give complete proofs/references below; general background references are~\cite{Borel-Tits65},\cite{Borel-Tits72},\cite{Borel-Tits73},\cite{Tits64},\cite{Margulis},\cite{Platonov-Rapinchuk}.

\medskip

\subsection{}
\label{sec_start}%
We fix temporarily an algebraic closure $\bbb F \supseteq F$; as is customary, we identify $\mathbf{GL}_d$ with the group of $\bbb F$\ti points $\mathbf{GL}_d(\bbb F)$. Consider the map $f:\Gamma\to\mathbf{GL}_d$ and recall that $\Gamma$ is finitely generated.

We claim that without loss of generality we may assume that we have a homomorphism $\Gamma\to\HH$ with infinite Zariski-dense image in some non-trivial connected adjoint absolutely simple group $\HH$. Indeed, upon replacing $\Gamma$ with a finite index subgroup (and the $G_i$'s accordingly), we can restrict ourselves to homomorphisms with Zariski-connected Zariski closure $\bbb{f(\Gamma)}$. Now we factor $\bbb{f(\Gamma)}$ by its radical, choose one of the almost simple factors of the quotient and pass to the associated adjoint group $\HH$. In order to prove the claim, it is enough to show that $\Gamma$ still had infinite image in one of the almost simple factors. If this is not the case, then $f(\Gamma)$ is virtually soluble (and still finitely generated); in particular, after possibly replacing $\Gamma$ with yet another finite index normal subgroup, and the $G_i$ accordingly, we would have a non-zero homomorphism $\Gamma\to\RR$. However, Shalom's superrigidity for characters~\cite[0.8]{Shalom00} implies that this homomorphism extends continuously to $G$; thus, some $G_i$ has a non-zero continuous homomorphism to $\RR$, which is incompatible with the definition of ``few factors''. This proves the claim.

\smallskip
\noindent
\textbf{Remark.~}The proof of~\cite[0.8]{Shalom00} is indeed carried out in the generality  of square-summable completely irreducible lattices in compactly generated locally compact groups.

\medskip

\subsection{}
\label{sec_K}%
Let $K$ be the subfield of $\bbb F$ generated by the traces in the adjoint representation of the images $\tau(\gamma)\in\HH$ of all elements $\gamma\in\Gamma$. This is a finitely generated field because $\Gamma$ is finitely generated and matrix multiplication is rational. Since we are in characteristic~$\neq 2,3$, there is no loss of generality in assuming that $\HH$ is defined over $K$ and that the image of $\Gamma$ lies in the $K$\ti points. Indeed, the assumption on the characteristic ensures that we are in the \emph{standard} case in the terminology of~\cite[VIII.3.15]{Margulis}. Therefore, the claim is a result of Vinberg, see~\cite[IX.1.8]{Margulis} (and~VIII.3.22 therein); we are thus reduced to study a homomorphism $\tau:\Gamma\to\HH(K)$ with infinite image.

\medskip

\subsection{}
Let $S$ be the set of all (inequivalent representatives of) valuations $v$ of $K$ such that the image of $\tau(\Gamma)$ is not relatively compact in $\HH(K_v)$ (for the $K_v$\ti topology); observe that this image is still Zariski-dense. Then $S\neq\varnothing$. The most forceful way to verify this is to apply a generalization of a lemma of Tits: Indeed, Breuillard-Gelander show~\cite[2.1]{Breuillard-Gelander} that any infinite subset of a finitely generated field becomes unbounded in a suitable completion.

\medskip

\subsection{}
\label{sec_isot_val}%
We claim that $S$ contains all the Archimedean valuations $v$ for which $\HH$ is $K_v$\ti isotropic (or equivalently, for which $\HH(K_v)$ is non-compact).

Indeed, suppose for a contradiction that $v$ is Archimedean, that $\HH(K_v)$ is non-compact but that the $v$\ti closure $C\se \HH(K_v)$ of $\tau(\Gamma)$ is compact. We distinguish two cases: Either $K_v \cong \RR$, and then Weyl's theorem on the algebraicity of compact subgroups of $\RR$\ti groups (see \emph{e.g.}~\cite[chap.~4, 2.1]{VinbergENC}) shows that the Zariski closure of $\tau(\Gamma)$ is $C\neq \HH(K_v)$, contradicting Zariski-density. The other possibility is $K_v\cong \CC$. Notice that $C$ is a Lie subgroup of $\HH(K_v)$; being compact, its traces in the adjoint representation of $\HH(K_v)$ are real~-- this follows \emph{e.g.} since $\mathrm{Ad}C$ preserves a real form of $\mathrm{Lie}(\HH(K_v))$, or alternatively because the Killing form is negative definite on $\mathrm{Lie}(C)$, see \emph{e.g.}~\cite[\S3]{Knapp97}. But then $K$ lies in $\RR$ under the isomorphism $K_v\cong \CC$, contradicting its density in $K_v$.

\medskip

\subsection{}
Let $v\in S$. The superrigidity theorem that we propose in~\cite{MonodCAT0} applies to $\Gamma \xrightarrow{\tau} \HH(K) \to \HH(K_v)$. More precisely, we recall for the reader's convenience:

\itshape
Let $\Gamma<G=G_1\times \cdots \times G_n$ be any square-summable weakly cocompact lattice in a product of locally compact $\sigma$\ti compact groups such that the projections of $\Gamma$ to each $G_i$ are dense. Let $\HH(K_v)$ be as before. Then any homomorphism $\Gamma\to\HH(K_v)$ with unbounded Zariski-dense image extends to a continuous homomorphism $G\to\HH(K_v)$ which moreover factors through some $G_i$.\upshape

(Apply Corollary~4 in~\cite{MonodCAT0}; Lemma~54 of~\cite{MonodCAT0} verifies that Zariski-density fulfills the assumptions of that corollary and Theorem~7 of~\cite{MonodCAT0} provides the extension to the non-cocompact case. Our current assumptions on the $G_i$'s are not needed in~\cite{MonodCAT0}, and neither is the part~(ii) of the definition of completely irreducible lattices.)

Therefore, there is an index $\delta(v)\in\{1,\ldots, n\}$ and a continuous homomorphism $\tau_v: G\to \HH(K_v)$ which factors through $G\twoheadrightarrow G_{\delta(v)}$ and extends $\tau$; we denote by $\sigma_v$ the resulting continuous homomorphism $\sigma_v:G_{\delta(v)} \to \HH(K_v)$. Observe that $\sigma_v$ is injective since otherwise it would have compact image by the assumption on the factor $G_{\delta(v)}$, contradicting the definition of $S$.

\medskip

\subsection{}
\label{sec_large_NA}%
Let $v\in S$. Since $G_{\delta(v)}$ has few factors, its identity component $G_{\delta(v)}^0$ is trivial or cocompact. In this subsection we assume the former and claim that $J\stackrel{\text{\rm\tiny def}}{=}\sigma_v(G_{\delta(v)})$ is a closed subgroup of $\HH(K_v)$ containing $\HH(K_v)^+$. In fact, it is enough to show that $J$ is open, since it is then closed and since any non-compact open subgroup contains $\HH(K_v)^+$ (the latter fact is stated in~\cite[9.10]{Borel-Tits73} and proved in~\cite{Prasad82}).

Since $G_{\delta(v)}$ is totally disconnected non-discrete (by Remark~\ref{sec_G_i}(c)), it admits a non-trivial open compact subgroup $C$, see~\cite[III \S4 \No~6]{BourbakiTGI}. Then $\sigma_v(C)$ is Zariski-dense because it is non-trivial and commensurated by the Zariski-dense subgroup $\tau(\Gamma)$. It is enough to show that $\sigma_v(C)$ is open, and we claim that this follows from the main result of~\cite{Pink98} (as in~\cite[p.~41]{Shalom00}). More precisely, it follows from~\cite[0.2]{Pink98} that there is a closed subfield $E<K_v$, an $E$\ti group $\HH_1$ and a $K_v$\ti isomorphism $\fhi:\HH\to\HH_1$ such that $\fhi\sigma_v(C)$ is an open subgroup of $\HH(E)$ (in fact this conclusion is already contained in the simpler statement~\cite[0.7]{Pink98} except for the discussion of the irreducibility of the adjoint representation~\cite[1.11]{Pink98}, which is relevant in positive characteristic). Since $\fhi(J)$ commensurates the Zariski-dense subgroup $\fhi\sigma_v(C)$, $\mathrm{Ad}\fhi(g)$ is defined over $E$ for all $g\in J$, and hence $\fhi(J)\se \HH_1(E)$ since $\HH_1$ is adjoint and $\mathrm{Ad}$ is defined over $E$~\cite[2.26]{Borel-Tits72}. Recalling that $\fhi$ preserves the traces~\cite[I.1.4.8]{Margulis}, we conclude $E=K_v$ from the definition of $K$. In conclusion, $\sigma_v(C)$ is open,  hence also $J$.

\medskip

\subsection{}
\label{sec_large_A}%
We consider now the case where $G_{\delta(v)}^0$ is cocompact and claim again that $\sigma_v(G_{\delta(v)})$ is closed and contains $\HH(K_v)^+$. Observe that $v$ is Archimedean since the image of $\sigma_v$ is non-compact.

By the solution of Hilbert's fifth problem there is a compact normal subgroup of $G_{\delta(v)}$ such that the quotient is isomorphic to a Lie group~\cite[4.6]{Montgomery-Zippin}; however in our case every compact normal subgroup is trivial. Thus $G_{\delta(v)}^0$ is a connected Lie group; we point out that it is semisimple by our assumption on $G_{\delta(v)}$. In any case, it is arcwise connected, and thus so is its image $L\stackrel{\text{\rm\tiny def}}{=}\sigma_v(G_{\delta(v)}^0)$. The Yamabe-Kuranishi theorem shows that $L$ is a Lie subgroup (see~\cite{Yamabe50}; for a short proof,~\cite{Goto69}). Since it is semisimple, it follows that $L$ is closed by a result of Mal'cev (see~\cite{Malcev45},\cite{Malcev46} or see~\cite[\S6]{Mostow50} for a detailed proof). As $G_{\delta(v)}^0$ is cocompact in $G_{\delta(v)}$, we deduce that $\sigma_v(G_{\delta(v)})$ is closed. Further, $L$ has finite index in its Zariski closure $\bbb L$, see~VIII.3.1-3 in~\cite{Hochschild81}. Since $\tau(\Gamma)$ normalizes $L$ and is Zariski-dense, it follows that $\HH(K_v)$ normalizes $\bbb L$. According to Tits~\cite{Tits64}, every (non-central) subgroup of $\HH(K_v)$ that is normalized by $\HH(K_v)^+$ contains $\HH(K_v)^+$ and thus $L\cap \HH(K_v)^+$ has finite index in $\HH(K_v)^+$, hence coincides with $\HH(K_v)^+$ by~\cite[6.7]{Borel-Tits73}. (This fact is more elementary here, since either $K_v = \RR$ and $\HH(K_v)^+$ is the identity component, or $K_v = \CC$ and $\HH(K_v)^+=\HH(K_v)$.) The claim follows.

Observe that by Baire's category theorem it follows from~\ref{sec_large_NA} and~\ref{sec_large_A} that $\sigma_v$ is a homeomorphism onto its image for all $v\in S$.

\medskip

\subsection{}
\label{sec_im_in}%
For any finite set $S'\se S$, there is a finite index normal subgroup $\Gamma_{S'}$ of $\Gamma$ whose image in $\HH(K_v)$ lies within $\HH(K_v)^+$ for all $v\in S'$. Indeed, since $S'$ is finite, it is enough to check this for each individual valuation $v$ for which $\HH$ is $K_v$\ti isotropic. But then, according to~\cite[6.14]{Borel-Tits73}, the group $\HH(K_v)/\HH(K_v)^+$, which is commutative, contains a finite index subgroup of finite exponent. Since $\Gamma$ is finitely generated, it follows that its image in $\HH(K_v)/\HH(K_v)^+$ is finite and thus we are done (we also used this argument in~\ref{sec_classical}).

Denote by $G_{\delta(v)}^{S'}$ the closure of the projection of $\Gamma_{S'}$ to $G_{\delta(v)}$; this closed subgroup has finite index. Since $\HH(K_v)^+$ is closed, we deduce further that $\sigma_v(G_{\delta(v)}^{S'})$ is in $\HH(K_v)^+$ for all $v\in S'$. Thus, by~\ref{sec_large_NA} and~\ref{sec_large_A}, and since $\HH(K_v)^+$ has no proper finite index closed subgroups~\cite{Tits64}, $\sigma_v(G_{\delta(v)}^{S'})=\HH(K_v)^+$ for all $v\in S'$.

\medskip

\subsection{}
The map $\delta:S\to \{1, \ldots, n\}$ is injective; in particular, $S$ is finite.

Indeed, assume that $v,w\in S$ have the same image $i$ and apply the discussion of~\ref{sec_im_in} to $S'=\{v,w\}$. Now both $\sigma_v$ and $\sigma_w$ induce isomorphisms of $G_i^{S'}$ onto their images $\HH(K_v)^+$ and respectively $\HH(K_w)^+$. We may now apply a result of Borel-Tits~\cite[8.13]{Borel-Tits73} (see also~I.1.8.2.III-IV in~\cite{Margulis}) to the end that the isomorphism $\sigma_w\sigma_v^{-1}: \HH(K_v)^+\to \HH(K_w)^+$ determines a topological isomorphism $\psi: K_v\to K_w$. However, since $\tau_v|_\Gamma = \tau = \tau_w|_\Gamma$ and thus $\sigma_w\sigma_v^{-1}$, which is defined on the whole of $\tau(\Gamma)$, preserves $\tau(\Gamma)$, it follows moreover that $\psi$ is the identity on $K$. In other words, $v$ and $w$ are equivalent, proving the claim.

\medskip

\subsection{}
\label{sec_global}%
$K$ is a global field. Indeed, it follows from its definition that it is infinite and finitely generated. If it is not global, then it has positive transcendence degree over its prime subfield, and in particular there is an element $\gamma\in\Gamma$ whose trace $\lambda\in K$ in the adjoint representation of $\HH$ is transcendental. Then, according to~\cite[IX.2.9]{Margulis}, $K$ admits infinitely many inequivalent valuations $v$ with $v(\lambda)$ arbitrarily large. However, if the subgroup generated by $\tau(\gamma)$ in $\HH(K_v)$ is relatively compact, then the absolute value of its trace is bounded by the dimension, so that all but finitely many of these valuations are in $S$ (compare~\cite[p.~306-307]{Margulis}). This contradicts the finiteness of $S$.

\medskip

\subsection{}
\label{sec_reduce_Gamma}%
Observe that thus far we have only once replaced $\Gamma$ with a finite index subgroup, and that was \emph{before} defining $K$ in~\ref{sec_K}. But with~\ref{sec_global} being now achieved, we shall not need to use the definition of $K$ in terms of $\Gamma$ anymore, so that we may now again freely replace $\Gamma$ with a finite index subgroup~-- and each $G_i$ accordingly.

In particular, since $S$ is finite, we may appeal to~\ref{sec_im_in} and replace $\Gamma$ with $\Gamma_S$ aswell as each $G_{\delta(v)}$ with $G_{\delta(v)}^S$ for all $v\in S$. Thus we may and shall henceforth assume that each $\sigma_v$ is a topological isomorphism $G_{\delta(v)}\to \HH(K_v)^+$.

\medskip

\subsection{}
\label{sec_finite_index}%
Let $H = \prod_{v\in S} \HH(K_v)$ and let $\Lambda<H$ be the image of $\HH(K(S))$ under the diagonal embedding $\Delta:\HH(K)\to H$. The claim of~\ref{sec_isot_val} is exactly the criterion from reduction theory ensuring that $\Lambda$ be a lattice in $H$ (see~\cite[I.3.2.5]{Margulis}). We claim that $\Delta\tau(\Gamma)\cap \Lambda$ has finite index in $\Delta\tau(\Gamma)$.

Recall first that for any non-Archimedean valuation $v$ of $K$, the subgroup $\HH(\O_v)$ of $\HH(K_v)$ is open, wherein $\O_v$ denotes the ring of integers of $K_v$. Therefore, if in addition $v\notin S$, then $\tau(\Gamma)\cap \HH(\O_v)$ has finite index in $\tau(\Gamma)$. It follows that for any finite set $\wt S\supseteq S$ of valuations,
$$\tau(\Gamma)\cap \HH(K(S))\ \text{ has finite index in }\ \tau(\Gamma)\cap \HH(K(\wt S)). \leqno{(\dagger)}$$
Since $\tau(\Gamma)$ is finitely generated and contained in $\HH(K)$, and since the latter is the union of the directed family $\HH(K(\wt S))$ where $\wt S$ ranges over the directed set of finite sets $\wt S \supseteq S$ of valuations, it follows that $\tau(\Gamma)$ is contained in $\HH(K(\wt S))$ for some such $\wt S$. Thus we deduce by~$(\dagger)$ that $\tau(\Gamma)\cap\HH(K(S))$ has finite index in $\tau(\Gamma)$, which is equivalent to the claim.

\medskip

\subsection{}
Since $\delta$ is injective, we can define a continuous homomorphism $\pi$ as the composition
$$\pi: G \twoheadrightarrow \prod_{v\in S} G_{\delta(v)} \xrightarrow{\ \prod \sigma_v} H^+ \stackrel{\text{\rm\tiny def}}{=}\prod_{v\in S} \HH(K_v)^+.$$
Then $\pi|_\Gamma = \Delta\tau$ and we claim that: (i)~$\pi(\Gamma)$ is commensurable with $\Lambda<H$, where (as in Remark~\ref{sec_thm_simple}(a)) we view $\pi(\Gamma)$ also as a subgroup of $H\supseteq H^+$; (ii)~$\delta$ is surjective onto $\{1, \ldots, n\}$, in particular $\pi$ is a topological isomorphism.

Indeed, $\Lambda$ being discrete in $H$, it follows that $\pi(\Gamma)\cap \Lambda$ is discrete in $H^+$; the latter having finite index in $\pi(\Gamma)$ by~\ref{sec_finite_index}, we see that $\pi(\Gamma)$ is discrete in $H^+$. Therefore, the second point in the definition of complete irreducibility~\ref{sec_G_i} shows that $\prod_{v\in S} G_{\delta(v)}$ cannot be a proper subproduct (since $\prod \sigma_v$ is a topological isomorphism), whence~(ii). Now it follows that $\pi(\Gamma)$ is a lattice in $H^+$; so is $\pi(\Gamma)\cap \Lambda$ by~\ref{sec_finite_index}. Since $H^+$ is normal cocompact in $H$, it has finite invariant covolume and thus any lattice in $H^+$ is a lattice in $H$ (see~\cite[1.6]{Raghunathan}). Thus $\pi(\Gamma)\cap \Lambda$ has also finite index in $\Lambda$, proving~(i).

The points~(i) and~(ii) just established complete the proof of Theorem~\ref{sec_few_factors} and hence of Theorem~\ref{sec_thm_simple} too.

\medskip

\subsection{}
\label{sec_cor}%
We turn to the corollary of~\ref{sec_thm_simple} and adopt its notation. We only need to show that one can reduce the problem to the case where $\Gamma$ is completely irreducible, since then Theorem~\ref{sec_thm_simple} applies. Thus, assume that one of the projections of $\Gamma$, say the projection to $G_1$, is discrete.

Consider the kernel $\Lambda$ of the (surjective) projection $\Gamma\to G_1$. Then $\Lambda$ is canonically realized as a closed subgroup of $G_2$; being normalized by the projection of $\Gamma$ to $G_2$, it is normal. We have now two cases in view of the topological simplicity of $G_2$: Either $\Lambda$ is trivial, and thus $\Gamma\cong G_1$ is simple. In that case, alternative~(ii) holds, since Mal'cev proved that every finitely generated linear group is residually finite~\cite{Malcev40}. The other case is $\Lambda=G_2$; then $G$ is discrete and hence $\Gamma=G_1\times G_2$ is a product of finitely generated simple groups and~(ii) holds as before by Mal'cev's result.

\medskip

\subsection{}
\label{sec_proof_few}%

Finally, we complete the picture of Section~\ref{sec_few_factors} by verifiying that every non-discrete locally compact group $G_i$ with property~$(*)$ has few factors. Observe that the groups $G_i$ considered in our theorems are indeed non-discrete because of the mere existence of a completely irreducible lattice, see Remark~\ref{sec_G_i}(c).

Let $G_i$ be a locally compact group with property~$(*)$. In particular, every non-trivial closed normal subgroup of $G_i$ is open. Therefore, the only way for the definition of ``few factors'' to fail is if there is a finite index open normal subgroup $H\lhd G_i$ and a non-zero continuous homomorphism $\chi:H\to\RR$. Since $\RR$ is torsion-free, property~$(*)$ implies that $\chi$ is injective and thus $H$ is a locally compact Abelian group without torsion. Since $H$ has no non-trivial closed normal subgroup with infinite quotient, it follows that $H$ is cyclic, hence countable. Now $G_i$ is also countable, hence discrete by Baire's theorem.

%
\def\cprime{$'$}

\end{document}